\documentclass[12pt]{article}
\usepackage[dvips]{graphicx}
\usepackage[latin1]{inputenc}
\usepackage{amssymb,amsmath,array}

\begin{document}

\noindent ESTIMATION AND TEST FOR MULTIDIMENSIONAL REGRESSION MODELS 
\vskip 3mm

\vskip 5mm
\noindent J. Rynkiewicz 

\noindent SAMOS-MATISSE

\noindent Centre d'Economie de la Sorbonne, 

\noindent Université de Paris I, Panthéon-Sorbonne, CNRS

\noindent 90 rue de Tolbiac 75013 Paris, France

\noindent rynkiewi@univ-paris1.fr

\vskip 3mm

\noindent ABSTRACT

This work is concerned with the estimation of multidimensional regression and the asymptotic behaviour of the test involved in selecting models. The main problem with such models is that we need to know the covariance matrix of the noise to get an optimal estimator. We show in this paper that if we choose to minimise  the logarithm of the determinant of the empirical error covariance matrix, then we get an asymptotically optimal estimator. Moreover, under suitable assumptions, we show that this cost function leads to a very simple asymptotic law for testing the number of parameters of an identifiable and regular regression model. Numerical experiments confirm the theoretical results. 

\vskip 3mm

\noindent 1. INTRODUCTION

Let us consider a sequence $\left(Y_{t},Z_{t}\right)_{t\in \mathbb{N}}$
of i.i.d. (i.e. independent identically distributed) random vectors. The law of $\left(Y_{t},Z_{t}\right)\in {\mathbb R}^{d}\times{\mathbb R}^{d'}$ is the same as the generic variable $(Y,Z)$.

We assume that the model can be written 
\[
Y_{t}=F_{w^0}(Z_{t})+\varepsilon _{t}, \eqno(1)
\]
where 

\begin{itemize}
\item $F_{w^0}$ is a parametric function, the true parameters being denoted $w^0$.
\item $(\varepsilon _{t})$ is an i.i.d.-centred noise with unknown definite positive
covariance matrix $\Gamma _{0}$.
\end{itemize}
The observations will be denoted with the lower-case letters $\left((z_t,y_t)\right)_{1\leq t\leq n}$.
This notation allows us to consider a wide range of regression models. For example $F_{w^0}$ can be a vector with $d$ lines and $d'$ columns, the parameter is the components of the matrix and the model will be a classical linear model. Another example can be constrained linear models, knowing that constraints can also be nonlinear. Finally, $F_{w^0}$ can also be a nonlinear parametric function like a multilayer perceptron (MLP). An MLP with $H$ hidden units (see Rumelhart et al. (1986)) is defined by a family of functions 
\[
F_w(z)=\sum_{h=1}^H b_js(a_j^Tz+c_j)+d,
\]
where $T$ denotes the transposition, $z\in {\mathbb R}^{d'}$, $s(t)=tanh(t)$ and \\
$w=\left(a_1,\cdots,a_H,b_1,\cdots,b_H,c_1,\cdots,c_H,d\right)\in {\mathbb R}^{Hd'}\times{\mathbb R}^{2H+1}$. We will focus on the MLP example, because it is a widely used tool for nonlinear regression (see White (1992)), but it could be any other non linear and differentiable model. 

 Note that, for an MLP function, there exists a finite number of transformations of the weights leaving these functions unchanged; these transformations form a finite group (see Sussmann (1992)). Therefore, we will consider equivalence classes of MLPs: two MLPs are in the same class if the first one is the image by such a transformation of the second one. The set of parameters considered is then the quotient space of parameters.  In the sequel, we will assume that the model is identifiable:
\[
F_{w_1}(Z)\stackrel{a.s.}{=}F_{w_2}(Z)\Leftrightarrow w_1=w_2. \eqno(2)
\]
 For example, this can be done if we consider one-hidden-layer MLPs with the true number of hidden units and with parameters in the quotient space. Another example is the linear regression function, with or without constraint. The consequence of the identifiability of the model is that, in most cases, the Hessian matrix of the model will be definite positive.Also, note that it is not hard to generalise all that is shown in this paper for stationary mixing variables and therefore for time series. For example, let us assume that the regression function verifies $\Vert F_{w^0}(z)\Vert\leq a\Vert z\Vert+b$, with $0\leq a<1$ and $b\in \mathbb R$.  Let $(Y_t)_{t\in \mathbb N}$ be the stationary solution  of the equation:
\[
Y_{t}=F_{w^0}(Y_{t-1})+\varepsilon _{t},
\]
where the noise $(\varepsilon _{t})_{t\in N}$ has a positive density everywhere with respect to the Lebesgue measure and is with an order moment strictly larger than 1. It is well known (see Duflo (1997)), that$(Y_t)_{t\in \mathbb N}$   will be geometrically ergodic and verifies a strong law of large numbers. In particular, as MLPs are bounded functions, if $F_{w^0}$ is an MLP function, all the proofs given in this paper will be valid exactly in the same way as in Yao (2000).

\vskip 3mm

\noindent 1.1. EFFICIENT ESTIMATION

The estimation of the model (1) is done by minimising a suitable cost function with respect to the parameters. 
A common choice for the cost function is the mean square error (MSE): 
\[
V_n(w):=\frac{1}{n}\sum _{t=1}^{n}\left\Vert y_{t}-F_{w}\left(z_{t}\right)\right\Vert ^{2},
\]
 where $\left\Vert .\right\Vert $ denotes the Euclidean norm on $\mathbb{R}^{d}$.
This function is widely used because in the linear case without constraint on the parameters this cost function is optimal (see L\"utkepohl (1993)). In fact, this cost function gives a satisfactory estimator when there is one and only one estimator which minimises the trace of the covariance matrix of the noise (see Magnus and Neudecker (1988)).
However in other cases (constraint linear model, non-linear regression, etc.) it leads in general to a suboptimal estimator (see, for example, L\"utkepohl (1993) for the constraint linear model).
Then, a better solution is to use an approximation of the covariance error matrix to compute the generalised least square estimator:
\[
\frac{1}{n}\sum _{t=1}^{n}\left(y_{t}-F_{w}\left(z_{t}\right)\right)^{T}\Gamma ^{-1}\left(y_{t}-F_{w}\left(z_{t}\right)\right),
\]

where $\Gamma $ has to be a good approximation of the true covariance matrix of the noise $\Gamma_{0}$. For example, if we use a sequence of matrices $\Gamma_n$ converging in probability to $\Gamma_0$, it is easy to show (see Chapter 5 in Galland (1987)) that the estimator obtained by minimising the cost function: 
\[
\frac{1}{n}\sum _{t=1}^{n}\left(y_{t}-F_{w}\left(z_{t}\right)\right)^{T}\Gamma_n^{-1}\left(y_{t}-F_{w}\left(z_{t}\right)\right)
\]
has the same asymptotic properties as the estimator which minimises 
\[
\frac{1}{n}\sum _{t=1}^{n}\left(y_{t}-F_{w}\left(z_{t}\right)\right)^{T}\Gamma_0^{-1}\left(y_{t}-F_{w}\left(z_{t}\right)\right).
\]
There are many ways to construct a sequence of $\left(\Gamma_{k}\right)_{k\in \mathbb{N}^{*}}$ yielding an approximation of $\Gamma _{0}$. The simplest is to use the ordinary least square estimator $\hat{W}^1_{n}:=\arg\min\frac{1}{n}\sum _{t=1}^{n}\left\Vert y_{t}-F_{w}\left(z_{t}\right)\right\Vert ^{2}$, in order to estimate the covariance matrix of the noise: \[
\Gamma^1_{n}:=\Gamma \left(\hat{W}^1_{n}\right):=\frac{1}{n}\sum _{t=1}^{n}(y_{t}-F_{\hat{W}^1_{n}}(z_{t}))(y_{t}-F_{\hat{W}^1_{n}}(z_{t}))^{T}.\]
Moreover, we can use this new covariance matrix to find a generalised
least square estimator $\hat{W}_{n}^{2}$: \[
\hat{W}_{n}^{2}=\arg \min _{w}\frac{1}{n}\sum _{t=1}^{n}\left(y_{t}-F_{w}\left(z_{t}\right)\right)^{T}\left(\Gamma^1_{n}\right)^{-1}\left(y_{t}-F_{w}\left(z_{t}\right)\right)\]
and calculate again a new covariance matrix \[
\Gamma^2_{n}:=\Gamma \left(\hat{W}_{n}^{2}\right)=\frac{1}{n}\sum _{t=1}^{n}(y_{t}-F_{\hat{W}_{n}^{2}}(z_{t}))(y_{t}-F_{\hat{W}_{n}^{2}}(z_{t}))^{T}.\]
Finally, it can be shown (see Gallant (1987)) that this procedure gives a sequence of parameters 
\[
\hat{W^1}_{n}\rightarrow \Gamma^1_{n}\rightarrow \hat{W}_{n}^{2}\rightarrow \Gamma^2_{n}\rightarrow \cdots 
\]
 minimising the logarithm of the determinant of the empirical covariance
matrix:
 
\[
U_{n}\left(w\right):=\log \det \left(\frac{1}{n}\sum _{t=1}^{n}(y_{t}-F_{w}(z_{t}))(y_{t}-F_{w}(z_{t}))^{T}\right).\eqno(3)
\]
Hence, the cost function $U_n(w)$ is the same as the generalised least square cost function with the best approximation of the true covariance matrix calculable with the available data.
Nevertheless it is important to note that the matrix $\Gamma_n$ is always a function of model parameters and it will be better to write $\Gamma_n(w)$ instead of $\Gamma_n$. The asymptotic study of the model must take into account the dependency of $\Gamma_n$ on these parameters, and the real function to study is in fact: 
\[
\frac{1}{n}\sum _{t=1}^{n}\left(y_{t}-F_{w}\left(z_{t}\right)\right)^{T}\Gamma_n^{-1}(w)\left(y_{t}-F_{w}\left(z_{t}\right)\right).
\]
This difficulty has always been overlooked except when the covariance matrix is included in the parameters of the model and this solution leads to consider a pseudo Gaussian likelihood as in Gourieroux et al. (1984). However, in this case it is necessary to reinforce the assumptions on the moment of the noise to obtain the asymptotic normality of the estimated covariance matrix. Although  the logarithm of the determinant of the empirical covariance is known to be related to the concentrated Gaussian likelihood function, it will be better to study it directly because such artificial strong assumptions about the noise are not needed.

For all these reasons, we propose to study the asymptotic properties of the cost function  $U_{n}\left(w\right)$ and  the estimator minimising this cost function: $\hat W_n :=\arg\min U_{n}\left(w\right)$ will be shown to have the same asymptotic behaviour as the generalised least square estimator using the true covariance matrix of the noise.

\vskip 3mm

\noindent 1.2. TESTING THE NUMBER OF PARAMETERS

The cost function $U_n(w)$ is not only optimal in the sense that it has the same asymptotic behaviour as the generalised least square estimator using the true covariance matrix of the noise, but it also leads to a very simple procedure for testing the nullity of the parameters.
Let $q$ be an integer smaller than $s$, we want to test ``$H_0: w\in \Theta_q \subset \mathbb R^q$'' against  ``$H_1: w\in \Theta_s \subset \mathbb R^s$'', where $\Theta_q$ and $\Theta_s$ are compact sets. $H_0$ expresses the fact that $w$ belongs to a subset of $\Theta_s$ with a lower parametric dimension than $s$ and so that $s-q$ parameters are equal to zero. 
If we consider the classic mean square error cost function: $V_n(w)$, we get the following test statistic (see Yao (2000)): 
\[
S_n=n\times\left(\min_{w\in \Theta_q}V_n(w)-\min_{w\in \Theta_s}V_n(w)\right).
\]
Under the null hypothesis $H_0$, it is shown in Yao (2000) that  $S_n$ converges in distribution to a weighted sum of $\chi^2_1$
\[
S_n\stackrel{\cal D}{\rightarrow}\sum_{i=1}^{s-q}\lambda_i\chi_{i,1}^2,
\]
where the $\chi_{i,1}^2$ are $s-q$ i.i.d. $\chi^2_1$ variables and $\lambda_i$ are strictly positive values, different from 1 if the true covariance matrix of the noise is not the identity matrix. 

However, if we use the function
\(
U_{n}\left(w\right)
\)
, under $H_0$, the test statistic:
\[
T_n=n\times\left(\min_{w\in \Theta_q}U_n(w)-\min_{w\in \Theta_s}U_n(w)\right)
\]
will converge to a classical $\chi^2_{s-q}$ and the asymptotic level of the test will be very easy to compute. This is another advantage of using the proposed cost function.
\paragraph{Organisation of the paper.}
In order to prove these properties, the paper is organised as follows: First, the main results are stated in three theorems, the first deals with the consistency of the estimator minimising $U_n(w)$, the second with its asymptotic normality and the third with the asymptotic law of the test procedure used to determine the number of parameters. Then, the theoretical results are confirmed by numerical experiments. The proofs of the theorems involving technical calculation of the first and second derivatives of $U_n(w)$ are postponed to the appendix.

\vskip 3mm

\noindent 2. ASYMPTOTIC PROPERTIES

First we give the conditions to state the consistency theorem, then  to state the asymptotic normality theorem of the estimator $\hat W_n$ minimising the function $U_n(w)$. In the sequel all the expectations will be calculated with respect to the true law of $(Y,Z)$.
\paragraph{Conditions for the consistency (C).}
\begin{enumerate}
\item The parameter space $\cal{W}$ is a compact space included in ${\mathbb R}^K$, with $K$ the dimension of the parameter vector $w$. The unique true parameter $w^0$ is assumed to be in the interior of  $\cal{W}$. Note that $w^0$ is unique because the model is assumed to be identifiable. The compactness of the parameter space means that parameters have to be bounded by a constant even if this constant can be very large. This is the case in practice if one uses a computer, since its numerical precision is finite. This rather classical hypothesis is needed to get the Glivenko-Cantelli property, which yields a uniform law of large numbers (see van der Vaart (1998)).   
\item The noise of the model $\varepsilon$ is square integrable. 
\item For almost all $z$ the function $w\mapsto F_w(z)$ is continuous; moreover  there exists a  square integrable function $m$ such that 
\[
\sup_{w\in \cal{W}}\Vert F_w(z) \Vert \leq m(z).
\] 
\end{enumerate}

These conditions are easily verifiable for regular models. For example, in the case of an MLP, it suffices to assume that the variable $Z$ is with finite third order moment. Indeed, in this case there exists a constant $C$ such that we have the following inequalities (see Yao (2000)):
\[
\begin{array}{l}
\sup_{w\in\cal{W}}\Vert F_w(Z) \Vert\leq C \\
\sup_{w\in\cal{W}}\Vert \frac{\partial F_w(Z)}{\partial w_k}\Vert\leq C(1+\Vert Z\Vert)\\
\sup_{w\in\cal{W}}\Vert \frac{\partial^2 F_w(Z)}{\partial w_k\partial w_l}\Vert\leq C(1+\Vert Z\Vert^2)\\
\sup_{w\in\cal{W}}\Vert \frac{\partial^3 F_{w}(Z)}{\partial w_j\partial w_k\partial w_l}\Vert\leq C \Vert(1+\Vert Z\Vert^3).
\end{array}
\]
For a linear model it suffices to assume that the variable $Z$ is with finite second order moment.

Then, we deduce  the theorem of consistency:
\paragraph{Theorem 1.}
Under the conditions (C), we have:
\[
\hat W_n\stackrel{a.s.}{\rightarrow}w^0.
\]

Now, we can establish the asymptotic normality for the estimator:

\paragraph{Conditions for the asymptotic normality (AN).}
\begin{enumerate}

\item There exists  a square integrable function $m_1$ such that, for all $k\in{1,\cdots,K}$:
\[
\sup_{w\in\cal{W}}\Vert \frac{\partial F_w(z)}{\partial w_k}\Vert\leq m_1(z).
\]

\item There exists integrable functions $m_2$ and $m_3$ such that for all $j,k,l\in{1,\cdots,K}$:
\[
\sup_{w\in\cal{W}}\Vert \frac{\partial^2 F_w(z)}{\partial w_j\partial w_k}\Vert\leq m_2(z)
\mbox{ and }
\sup_{w\in\cal{W}}\Vert \frac{\partial^3 F_w(z)}{\partial w_j\partial w_k\partial w_l}\Vert\leq m_3(z).
\]
\end{enumerate}

Thus, we deduce the theorem: 
\paragraph{Theorem 2.}
Under the conditions (C) and (AN), when $n\rightarrow\infty$,
\[
\sqrt{n}(\hat W_n-w^0)\stackrel{\cal D}{\rightarrow}{\cal N}(0,I_0^{-1}),
\]
where, if we note
\[
B(w_k,w_l):=\frac{\partial F_{w}(Z)}{\partial w_{k}}\frac{\partial F_{w}(Z)}{\partial w_{l}}^T,
\]
the component $(k,l)$ of the matrix $I_0$ is:
\[
tr\left(\Gamma^{-1}_0E\left(B(w^0_k,w^0_l)\right) \right).
\]

\paragraph{Remark.} 
If $W_n^*$ is the estimator of the generalised least squares:
\[
W_n^*:=\arg\min \frac{1}{n}\sum _{t=1}^{n}\left(Y_{t}-F_{w}\left(Z_{t}\right)\right)^{T}\Gamma_0^{-1}\left(Y_{t}-F_{w}\left(Z_{t}\right)\right),
\]
then it is easy to check that  
\[
\sqrt{n}(W_n^*-w^0)\stackrel{\cal D}{\rightarrow}{\cal N}(0,I_0^{-1}).
\]

So, $\hat W_n$ has the same asymptotic behaviour as the generalised least square estimator with the true covariance matrix  $\Gamma^{-1}_{0}$ which is asymptotically optimal (see for example Ljung (1999)).  Therefore, the proposed estimator is asymptotically optimal too.
\paragraph{Asymptotic distribution of the test statistic $T_n$.}
Let us remind that we want to test ``$H_0: w\in \Theta_q \subset \mathbb R^q$'' against  ``$H_1: w\in \Theta_s \subset \mathbb R^s$''. $H_0$ expresses the fact that $w$ belongs to a subset of $\Theta_s$ with a parametric dimension smaller than $s$ so that $s-q$ parameters are equal to zero. 

Let us  write\\  $\hat W_n=\arg\min_{w\in \Theta_s}U_n(w)$ and
$\hat W^0_n=\arg\min_{w\in \Theta_q} U_n(w)$, where $\Theta_q$ is viewed as a subset of $\Theta_s$. Under the null hypothesis $H_0$, the asymptotic distribution of $T_n$ is a consequence of Theorem 2. Indeed, if we replace $T_n$ by its Taylor expansion around $\hat W_n$ and $\hat W^0_n$, following  van der Vaart (1998), chapter 16, we have the classical development:
\paragraph{Theorem 3.}
under the conditions (C) and (AN), if we assume that matrix $I_0$ is not singular and under the null hypothesis $H_0$, we have:
\[
T_n=n\left(\hat W_n-\hat W^0_n\right)^T I_0\left(\hat W_n-\hat W^0_n\right)+o_P(1)\stackrel{\cal D}{\rightarrow}\chi^2_{s-q},
\]
where $o_P(1)$ means ``negligible in probability'' and is defined in van der Vaart (1998).
\vskip 3mm

\noindent 3. EXPERIMENTAL RESULTS

\vskip 3mm

\noindent 3.1. SIMULATED EXAMPLE

Although the estimator associated with the cost function $U_{n}\left(w\right)$
is theoretically better than the ordinary mean least square estimator,
we have to confirm this fact by simulation in some cases. For example, there are some pitfalls in practical situations with MLPs. 

The first point is that we have no guarantee of reaching the global minimum of the cost function because we use differential optimisation to estimate $\hat W_n$. Hence, we can only hope to find a good local minimum if we use many estimations with different initial weights.

The second point is the fact that MLPs are black boxes, which means that it is difficult to interpret on their parameters and it is almost impossible to compare MLPs by comparing their parameters, even if we try to take into account the possible permutations of the weights. 

These reasons explain why we choose to compare the estimated covariance
matrices of the noise instead of directly comparing the estimated parameters of MLPs.

\paragraph{The model.}

To simulate our data, we use an MLP with 2 inputs, 3 hidden units,
and 2 outputs. We choose to simulate a time series, because it is a very
easy task as the outputs at time $t$ are the inputs
for time $t+1$. Moreover, with MLPs, the statistical properties
of such a model are the same as with independent identically distributed
(i.i.d.) data. Indeed, since the MLP function is bounded and the noise has a density positive everywhere with respect to the Lebesgue measure,  the time series simulated is an example of a process with a geometrically ergodic solution (see Yao (2000)) and verifies a strong law of large numbers.

The equation of the model is the following 
\[
Y_{t+1}=F_{w_{0}}\left(Y_{t}\right)+\varepsilon _{t+1},
\]

where 

\begin{itemize}
\item $Y_0=(0,0)$.
\item $\left(Y_{t}\right)_{1\leq t\leq 1000}$, $Y_{t}\in \mathbb{R}^{2}$,
is the bi-dimensional simulated random process. 
\item $F_{w_{0}}$ is an MLP function with weights $w_{0}$ randomly chosen between $-2$ and 2.
\item $(\varepsilon _{t})$ is an i.i.d. Gaussian centred noise with covariance matrix $\Gamma _{0}=\left(\begin{array}{cc}
 1.81 & 1.8\\
 1.8 & 1.81\end{array}
\right)$.
\end{itemize}
In order to empirically study the statistical properties of our estimator,
we make $100$ independent simulations of the bi-dimensional time series
of length $1000$. 
\paragraph{Results.}
For each time series we estimate the weights of the MLP using the cost
function $U_n\left(w\right)$ and the ordinary least square estimator. 
The estimations were made using the second order algorithm BFGS (see Press et al. (1992)), and for each estimation we chose the best result obtained
after $20$ random initialisations of the weights in the hope of avoiding 
 to plague our learning with poor local minima. 

We here show the mean of the estimated covariance matrices of the noise for
$U_n(w)$ and the mean square error (MSE) cost function: \[
U_n\left(w\right)\,:\, \left(\begin{array}{cc}
 1.793 & 1.785\\
 1.785 & 1.797\end{array}
\right)\mbox {\, and\, MSE}\,:\, \left(\begin{array}{cc}
 1.779 & 1.767\\
 1.767 & 1.783\end{array}
\right).\]
The estimated standard deviation of the terms of the matrices are
all equal to $0.003$, so the differences observed between the terms of the two matrices are greater than twice their standard deviation and probably not due to chance. 
We can see that the estimated covariance of the noise is on average better with
the estimator associated with the cost function $U_n\left(w\right)$.
In particular, it seems that there is slightly less over-fitting with
this estimator, and the non-diagonal terms are greater than with the least square estimator. As expected, the determinant of the mean matrix associated with $U_n(w)$ is 0.036 instead of 0.050 for the matrix associated with the MSE. 
\vskip 3mm

\noindent 3.2. APPLICATION TO REAL TIME SERIES: POLLUTION OF OZONE

Ozone is a reactive oxide, which is formed both in the stratosphere and troposphere. Near the surface of the ground, ozone is directly harmful to human health, plant life and damages physical materials. The population, especially in large cities and in suburban zones which suffer from summer smog, wants to be warned of high pollutant concentrations in advance. Statistical ozone modelling and more particularly regression models have been widely studied, see Comrie (1997), Gardner and Dorling (1998). Generally, linear models do not seem to capture all the complexity of this phenomenon. Thus, the use of nonlinear techniques is recommended to deal with ozone prediction. Here we want to predict ozone pollution at two sites at the same time. The sites are in the south of Paris (13th district) and at the top of the Eiffel Tower. As these sites are very near each other we can expect that the two components of the noise are highly correlated. 
\paragraph{The model.}
The neural model used in this study is autoregressive and includes exogenous parameters (called NARX model). 
Our aim is to predict the maximum level of ozone pollution of the next day, given today's maximum level of pollution and the maximal temperature of the next day. If we note $Y^1$ the maximum level of pollution for Paris 13, $Y^2$ the maximum level of pollution for the Eiffel Tower and $Temp$ the temperature, the model can be written as follows: 
\[
(Y^1_{t+1},Y^2_{t+1})=F_w(Y^1_t,Y^2_t,Temp_{t+1})+\varepsilon_{t+1}.\eqno(4)
\]
If we assume that the temperature $(Temp_{t})_{t\in \mathbb N}$ is a geometrically ergodic process and that the noise $(\varepsilon_t)_t\in\mathbb N$ has a strictly positive density with respect to the Lebesgue measure, as $F_w$ is a bounded function, the stationary solution $(Y^1_{t},Y^2_{t})_{t\in \mathbb N}$  of the equation (4) will be geometrically ergodic and the previous theorems can easily be extended to this time series.
  
As usual with real time series, over-training is a crucial problem. MLPs are very over-parametrised models. This occurs when the model learns the details of the noise of the training data. Over-trained models have very poor performance on fresh data. In this study, to avoid over-training, we use the Statistical Stepwise Method (SSM) pruning technique, using a BIC-like information criterion (Cottrell et al. (1995)). The MLP with the minimal dimension is found by eliminating of the irrelevant weights in order to minimise a BIC-like criterion, that is to say the cost function penalised by the term $q\frac{\ln(n)}{n}$, where $q$ is the number of parameters of the model and $n$ is the number of observations. Here, we will compare the behaviour of this method for both cost functions: the mean square error (MSE) and the logarithm of the determinant of the empirical covariance matrix of the noise ($U_n(w)$). 
\paragraph{Dataset.}
This study presents the ozone concentration of the Air Quality Network of the Ile de France Region (AIRPARIF, Paris, France). The data used in this work span from 1994 to 1997.  According to the model, we have the following explicative variables:
\begin{itemize}
\item The maximum temperature of the day
\item Persistence is used by introducing the previous day's ozone peak.
\end{itemize}
Before being used in the neural network, all these data have been centred and normalised. 
The learning dataset consists of observations from 1994 to 1996. Only the months from April to September are used because there is no peak of ozone pollution during the winter period.   The months from April to September of 1997 are kept for a test set, which will be used for evaluating models.
\paragraph{The results.}
For the learning set, we get the following results:
\[
U_n(w)\,:\, \left(\begin{array}{cc}
 0.26 & 0.19\\
 0.19 & 0.34\end{array}
\right)\mbox {\, and\, MSE}\,:\, \left(\begin{array}{cc}
 0.26 & 0.18\\
 0.18 & 0.34\end{array}
\right).\]
For the test set, we get the following results:
\[
U_n(w)\,:\, \left(\begin{array}{cc}
 0.32 & 0.21\\
 0.21 & 0.39\end{array}
\right)\mbox {\, and\, MSE}\,:\, \left(\begin{array}{cc}
 0.34 & 0.20\\
 0.20 &0.41\end{array}
\right).\]
The two matrices are almost the same for the learning set,  however the non-diagonal terms are greater for the $U_n(w)$ cost function.
The best MLP for $U_n(w)$ has 13 weights, and the best MLP for the MSE cost function has 15 weights. Hence, the proposed cost function leads to a more parsimonious model, certainly because the pruning technique is very sensitive to the variance of estimated parameters. This gain is valuable regarding the generalisation capacity of the model, since in this way the difference is almost null for the learning data set but is greater for the test data. For comparison, we did the training with a one-output MLP to predict each level of pollution and the results match the diagonal terms of the MSE cost function.

\vskip 3mm

\noindent 4. CONCLUSION

In the linear multidimensional regression model without constraint the optimal estimator has an analytic solution and it minimises both the ordinary mean square function and $U_n(w)$, therefore it is not useful, for this case, to consider $U_n(w)$. However, for the constrained linear model and for the non-linear multidimensional regression model, the ordinary least square estimator is sub-optimal if the covariance matrix of the noise is not the
identity matrix. We can overcome this difficulty by using the cost function $U_n(w)=\log\det(\Gamma_n(w))$. Indeed, this cost function is the same as the generalised least square cost function with the best approximation of the true covariance matrix calculable with the available data. In this paper, the proofs of  the consistency, of the normality asymptotic and of the optimality of this estimator have been provided. Moreover, we have proved that, if the model is identifiable, this cost function leads to a simpler test to determine the number of weights. These theoretical results have been confirmed by a simulated example, and we have seen on a real time series that we can expect slight improvement, especially in model selection, because pruning techniques are very sensitive to the variance of the estimated weights. Nevertheless, we have to note that the main difficulty in regression with MLPs is the lack of theoretical justification for procedures determining the number of hidden units. Indeed, determining the true number of hidden units is very important in order to have an identifiable model. For practical situations and without theoretical justification, a BIC-like penalised cost function seems to work well. 

\vskip 3mm

\noindent 5. APPENDIX

\paragraph{Proof of theorem 1.} 
First we have to show that the limit, as $n$ goes to infinite, of $U_n(w)$ is minimised only by $w^0$. 
\paragraph{Lemma 1.}
Under the conditions (C):
\[
\lim_{n\rightarrow \infty}U_n(w)-U_n(w^0)\stackrel{a.s.}{\geq}0
\]
and 
\[
\lim_{n\rightarrow \infty}U_n(w)-U_n(w^0)\stackrel{a.s.}{=}0\Leftrightarrow w=w^0.
\]
\paragraph{Proof:}
Let us note 
\[
\Gamma(w)=E\left((Y-F_w(Z))(Y-F_w(Z))^T\right) \eqno(5)
\]
the expectation of the covariance matrix of the noise for the model parameter $w$ and remark that $\Gamma_0=\Gamma(w^0)$.
By the strong law of large numbers we have
\[
\begin{array}{l}
U_n(w)-U_n(w^0)\stackrel{a.s.}{\rightarrow}\log \det(\Gamma(w))-\log \det(\Gamma_0)=\log\frac{\det(\Gamma(w))}{\det(\Gamma_0)}\\
=\log \det\left(\Gamma^{-1}(w^0) \left(\Gamma(w)-\Gamma_0\right)+I_d\right),
\end{array}
\]
where $I_d$ denotes the identity matrix of $\mathbb R^d$.
So,  the lemma is true if $\Gamma(w)-\Gamma_0$ is a positive matrix, null only if $w=w^0$, because the determinant of $\left(\Gamma^{-1}(w^0) \left(\Gamma(w)-\Gamma_0\right)+I_d\right)$ will be bigger than 1. This is the case since
\[
\begin{array}{l}
\Gamma(w)=E\left((Y-F_w(Z))(Y-F_w(Z))^T\right)\\
=E\left((Y-F_{w^0}(Z)+F_{w^0}(Z)-F_w(Z))(Y-F_{w^0}(Z)+F_{w^0}(Z)-F_w(Z))^T\right)\\
=E\left((Y-F_{w^0}(Z))(Y-F_{w^0}(Z))^T\right)+E\left((F_{w^0}(Z)-F_w(Z))(F_{w^0}(Z)-F_w(Z))^T\right)\\
=\Gamma_0+E\left((F_{w^0}(Z)-F_w(Z))(F_{w^0}(Z)-F_w(Z))^T\right).
\end{array}
\]
Then, the lemma is proved because the model is assumed to be identifiable (see equation (2)), so
\[
E\left((F_{w^0}(Z)-F_w(Z))(F_{w^0}(Z)-F_w(Z))^T\right),
\]
is a positive matrix, null only if $w=w^0$ $\blacksquare$

From the assumptions (C), following example 19.8 of van der Vaart (1998) the set of functions
\[
w\mapsto(Y-F_w(Z))(Y-F_w(Z))^T,\ w\in{\cal W}
\] is Glivenko-Cantelli.

Now, by lemma 1, we remark that for all neighbourhood $\cal O$ of $w^0$ there exists a number $\eta({\cal O})>0$ such that for all $w\notin \cal O$ we have
\[
\log\det\left(\Gamma(w)\right)>\log\det\left(\Gamma_0\right)+\eta({\cal O}).
\]

In order to show the strong consistency we have to prove that for all neighbourhood $\cal O$ of $w^0$ we have 
\(
\lim_{n\rightarrow \infty}\hat W_n \stackrel{a.s.}{\subset}{\cal O}
\), which is equivalent to 
\[
\lim_{n\rightarrow \infty}\log\det\left(\Gamma(\hat W_n)\right)-\log\det\left(\Gamma_0\right)\stackrel{a.s.}{<}\eta({\cal O}),
\]
where $\Gamma(\hat W_n)$ is defined by equation (5).
By definition, we have:
\[
\log\det\left(\Gamma_n(\hat W_n)\right)\stackrel{a.s.}{\leq} \log\det\left(\Gamma_n(w^0)\right).
\]
The Glivenko-Cantelli  property and the continuity of the function $\Gamma\mapsto \log\det(\Gamma)$ imply that
\[ 
\lim_{n\rightarrow \infty} \log\det\left(\Gamma_n(w^0)\right)-\log\det\left(\Gamma_0\right)\stackrel{a.s.}{=}0,
\]
therefore
\[
\lim_{n\rightarrow \infty} \log\det\left(\Gamma_n(\hat W_n)\right)\stackrel{a.s.}{<}\log\det\left(\Gamma_0\right)+\frac{\eta({\cal O})}{2}.
\]
They also imply that
\[ 
\lim_{n\rightarrow \infty} \log\det\left(\Gamma_n(\hat W_n)\right)-\lim_{n\rightarrow \infty}\log\det\left(\Gamma(\hat W_n)\right)\stackrel{a.s.}{=}0
\]
and finally
\[
\lim_{n\rightarrow \infty} \log\det\left(\Gamma(\hat W_n)\right)-\frac{\eta({\cal O})}{2}\stackrel{a.s.}{<}
\lim_{n\rightarrow \infty}\log\det\left(\Gamma_n(\hat W_n)\right)\stackrel{a.s.}{<}\log\det\left(\Gamma_0\right)+\frac{\eta({\cal O})}{2},
\]
so
\[
\lim_{n\rightarrow \infty} \log\det\left(\Gamma(\hat W_n)\right)\stackrel{a.s.}{<}\log\det\left(\Gamma_0\right)+\eta({\cal O})
\]
$\blacksquare$

\paragraph{Proof of theorem 2.} 
As usual, the asymptotic normality of the estimator minimising $U_n(w)$ is a consequence of the Taylor expansion of $U_n(w)$ around the parameter $w^0$. So, the computation of the first and second derivative of $U_n(w)$ is necessary to get these results.

Let us introduce a notation:  if $F_w(z)$ is a $d$-dimensional parametric function depending on a parameter vector $w$, we write $\frac{\partial F_w(z)}{\partial w_k}$ (resp. $\frac{\partial^2 F_w(z)}{\partial w_k\partial w_l}$) for the $d$-dimensional vector of the partial derivative (resp. second order partial derivatives) of each component of $F_w(z)$. Moreover, if $\Gamma(w)$ is a matrix depending on $w$, let us write $\frac{\partial}{\partial w_k}\Gamma(w)$ the matrix of partial derivatives of each component of $\Gamma(w)$. 
\subparagraph{First derivatives.}
Now, if $\Gamma_n(w)$ is a matrix depending on the parameter vector $w$, we get  from Magnus and Neudecker (1988)
\[
\frac{\partial }{\partial w_{k}}\log \det \left(\Gamma _{n}(w)\right)=tr\left(\Gamma_{n}^{-1}(w) \frac{\partial }{\partial w_{k}}\Gamma _{n}(w)\right),
\]
with 
\[
\Gamma_n(w)=\frac{1}{n}\sum _{t=1}^{n}(y_{t}-F_{w}(z_{t}))(y_{t}-F_{w}(z_{t}))^{T}.
\] 
Note that this matrix $\Gamma_{n}(w)$ and its inverse are symmetrical. 
Now, if we write 
\[
A_n(w_k)=\frac{1}{n}\sum _{t=1}^{n}\left(-\frac{\partial F_{w}(z_{t})}{\partial w_{k}}(y_t-F_{w}(z_{t}))^T\right),
\eqno(6)
\]
using the fact that 
\[
tr\left(\Gamma_{n}^{-1}(w)A_n(w_k)\right)=tr\left(A_n^T(w_k)\Gamma_{n}^{-1}(w)\right)=tr\left(\Gamma_{n}^{-1}(w)A_n^T(w_k)\right),
\]
we get 
\[
\frac{\partial }{\partial w_{k}}\log \det \left(\Gamma _{n}(w)\right)=2tr\left(\Gamma_{n}^{-1}(w)A_n(w_k)\right).
\]
As we will see in an example, the calculation of this derivative is generally easy.
\subparagraph{Example: calculation of the derivative for an MLP.}
The $i$th component of a multidimensional function will be denoted $F_{w}(z_{t})(i)$ and for a matrix $A=\left(A_{ij}\right)$, we write $vec(A)$ the vector obtained by concatenation of the columns of $A$. Following the previous results, and according to Magnus and Neudecker (1988), we can write the derivative of $\log (\det (\Gamma _{n}\left(w\right)))$
with respect to the weight $w_{k}$: \[
\frac{\partial }{\partial w_{k}}\log \det (\Gamma _{n}\left(w\right))=vec\left(\Gamma^{-1}_{n}\left(w\right)\right)^{T}vec\left(\frac{\Gamma _{n}\left(w\right)}{\partial w_{k}}\right),\eqno(7)\]
 with $\frac{\Gamma _{n}\left(w\right)}{\partial w_{k}}$ the matrix whose component $ij$ is:
\[
\frac{1}{n}\sum _{t=1}^{n}\left[-\frac{\partial F_{w}(z_{t})(i)}{\partial w_{k}}\times \left(y_{t}-F_{w}(z_{t})\right)(j)-\frac{\partial F_{w}(z_{t})(j)}{\partial w_{k}}\left(y_{t}-F_{w}(z_{t})\right)(i)\right].\eqno(8)
\]

Back-propagation is the standard way to compute the derivatives with an MLP function (see Haykin (1999)).
Here,  if we consider the MLP restricted to the output $i$, the quantity $\frac{\partial F_{w}(z_{t})(i)}{\partial w_{k}}$ can be computed by back-propagating the constant $1$. For example, figure \ref{cap:MLP-restricted-to} gives an example of an MLP restricted to the output $2$.

\begin{figure}[h]

\caption{\label{cap:MLP-restricted-to}MLP restricted to the output $2$:
the continuous lines}

HERE, Figure 1

\end{figure}

Hence, the computation of the gradient of $U_{n}\left(w\right)$ with respect
to the parameters of the MLP is straightforward. We have to compute
the derivative with respect to the weights of each single output MLP extracted from the original MLP.
This computation can be done by back-propagating the constant value $1$. 
Then, according to formula (8), we compute the derivative
of each term of the empirical covariance matrix of the noise. Finally the gradient
is obtained by the sum of all the derivative terms of the empirical
covariance matrix multiplied by the terms of its inverse as in formula (7).
\subparagraph{Second derivatives.}
We  now write
\[
B_n(w_k,w_l):=\frac{1}{n}\sum _{t=1}^{n}\left( \frac{\partial F_{w}(z_{t})}{\partial w_{k}}\frac{\partial F_{w}(z_{t})}{\partial w_{l}}^T\right)
\]
and
\[
C_n(w_k,w_l):=\frac{1}{n}\sum _{t=1}^{n}\left( -(y_t-F_{w}(z_{t})) \frac{\partial^2 F_{w}(z_{t})}{\partial w_{k}\partial w_{l}}^T\right).
\]
We get

\[
\begin{array}{l}
\frac{\partial^2 U_n(w)}{\partial w_k\partial w_l}=\frac{\partial }{\partial w_{l}}2tr\left(\Gamma_{n}^{-1}(w)A_n(w_k)\right)=\\
2tr\left(\frac{\partial\Gamma_{n}^{-1}(w)}{\partial w_{l}}A(w_k)\right)+2tr\left(\Gamma_{n}^{-1}(w)B_n(w_k,w_l)\right)+2tr\left(\Gamma_{n}(w)^{-1}C_n(w_k,w_l)\right).\\
\end{array}
\]
Now, Magnus and Neudecker (1988) give an analytic form of the derivative of an inverse matrix: 
\[
\frac{\partial \Gamma^{-1}(w)}{\partial w_k}=-\Gamma^{-1}(w)\left(\frac{\partial \Gamma(w)}{\partial w_k}\right)\Gamma^{-1}(w)
\]
so
\[
\begin{array}{l}
\frac{\partial^2 U_n(w)}{\partial w_k\partial w_l}=2tr\left(\Gamma_{n}^{-1}(w)\left(A_n(w_k)+A_n^T(w_k)\right)\Gamma_{n}^{-1}(w)A_n(w_k)\right)+\\
2tr\left(\Gamma_{n}^{-1}(w)B_n(w_k,w_l)\right)+2tr\left(\Gamma_{n}^{-1}(w)C_n(w_k,w_l)\right)
\end{array}
\eqno(9)
\]

Now, theorem 2 will follow from this fundamental lemma.
\paragraph{Lemma 2.}

Let $\Delta U_n(w^0)$ be the gradient vector of $U_n(w)$ at $w^0$ and $HU_n(w^0)$ be the Hessian matrix of  $U_n(w)$ at $w^0$. 

We finally define
\[
B(w_k,w_l):=\frac{\partial F_{w}(Z)}{\partial w_{k}}\frac{\partial F_{w}(Z)}{\partial w_{l}}^T.
\]

Then, under the assumption (AN) we get:  
\begin{enumerate}
\item $\sqrt{n}\Delta U_n(w^0)\stackrel{\cal D}{\rightarrow}{\cal N}(0,4I_0)$
\item $HU_n(w^0)\stackrel{a.s.}{\rightarrow}2I_0$,
\end{enumerate}
where the component $(k,l)$ of the matrix $I_0$ is:
\[
tr\left(\Gamma^{-1}_0E\left(B(w^0_k,w^0_l)\right) \right).
\]

\paragraph{Proof of lemma 2:}
Let us begin with the first result.
Under the condition (AN)-1, $A_n(w^0_k)$ (see equation (6)) is square integrable, so it verifies the central limit theorem. As the $k$th term of $\Delta U_n(w^0)$ is equal to 
\(2tr\left(\Gamma_{n}^{-1}(w^0)A_n(w^0_k)\right)\)
and $\Gamma_{n}^{-1}(w^0)\stackrel{a.s.}{\rightarrow}\Gamma_0^{-1}$,  by the Slutsky lemma (see lemma 2.8 of van der Vaart (1998)), $\Delta U_n(w^0)$ verifies the central limit theorem too.  Now,  let us write
\[
A(w_k)=\left(-\frac{\partial F_{w}(Z)}{\partial w_{k}}(Y-F_{w}(Z))^T\right)
\] 
and  write $\frac{\partial U(w)}{\partial w_k}:=\frac{\log\det(\Gamma(w))}{\partial w_k}$ .
We first remark that the component $(k,l)$ of the matrix $4I_0$ is:
\[
E\left(\frac{\partial U(w^0)}{\partial w_k}\frac{\partial U(w^0)}{\partial w^0_l}\right)=E\left(2tr\left(\Gamma^{-1}_0A^T(w^0_k)\right)\times2tr\left(\Gamma^{-1}_0A(w^0_l)\right)\right)
\]
and, since the trace of the product is invariant by circular permutation, 
\[
\begin{array}{l}
E\left(\frac{\partial U(w^0)}{\partial w_k}\frac{\partial U(w^0)}{\partial w^0_l}\right)\\
=4E\left( -\frac{\partial F_{w^0}(Z)^T}{\partial w_k}\Gamma^{-1}_0(Y-F_{w^0}(Z))(Y-F_{w^0}(Z))^T\Gamma^{-1}_0\left(-\frac{\partial F_{w^0}(Z))}{\partial w_l}\right)\right)\\
=4E\left(\frac{\partial F_{w^0}(Z)^T}{\partial w_k}\Gamma^{-1}_0\frac{\partial F_{w^0}(Z)}{\partial w_l}\right)\\
=4tr\left(\Gamma^{-1}_0E\left(\frac{\partial F_{w^0}(Z)}{\partial w_k}\frac{\partial F_{w^0}(Z)^T}{\partial w_l}\right) \right)\\
=4tr\left(\Gamma^{-1}_0E\left(B(w^0_k,w^0_l)\right)\right). 

\end{array}
\]
This proves the first result.

Let us now prove the second result.  For the component $(k,l)$ of the expectation of the Hessian matrix, we remark that 
\[
\lim_{n\rightarrow \infty}A_n(w^0_k)=E\left(A(w^0_k)\right)=0
\]
because the noise $\varepsilon=Y-F_{w^0}(Z)$ is centred and independent of the random variable $Z$.
Hence
\[
\begin{array}{l}
\lim_{n\rightarrow \infty}tr\left(\Gamma_{n}^{-1}(w^0)A_n(w^0_k)\Gamma_{n}^{-1}(w^0)A_n(w^0_k)\right)=\\
\lim_{n\rightarrow \infty}tr\left(\Gamma_{n}^{-1}(w^0)A^T_n(w^0_k)\Gamma_{n}^{-1}(w^0)A_n(w^0_k)\right)=0
\end{array}
\]
and, for the same reason
\[
\lim_{n\rightarrow \infty}tr\Gamma_{n}^{-1}C_n(w^0_k,w^0_l)=0.
\]
Finally
\[
\begin{array}{l}
\lim_{n\rightarrow \infty}HU_n(w^0)=\lim_{n\rightarrow \infty}2tr\left(\Gamma_{n}^{-1}(w^0)\left(A_n(w^0_k)+A^T_n(w^0_k)\right)\Gamma_{n}^{-1}(w^0)A_n(w^0_k)\right)+\\
2tr\Gamma_{n}^{-1}(w^0)B_n(w^0_k,w^0_l)+2tr\Gamma_{n}^{-1}C_n(w^0_k,w^0_l)\\
=2tr\left(\Gamma^{-1}_0E\left(B(w^0_k,w^0_l)\right)\right)\\
\blacksquare
\end{array}
\]
As the matrix $I_0$ is assumed to be invertible, following the same argument of local asymptotic normality as in Yao (2000), we get the Taylor formula with an integral remainder:
\[
\Delta U_n(\hat W_n)=\Delta U_n(w^0)+\int_0^1HU_n(\hat W_n+u(\hat W_n-w^0)du.
\]  
The condition (AN)-2 implies that
\[
\Vert \frac{\partial F_{w_1}(z)}{\partial w_k}-\frac{\partial F_{w_2}(z)}{\partial w_k}\Vert\leq \Vert w_1 -w_2\Vert m_2(z)
\]
and
\[
\Vert \frac{\partial^2 F_{w_1}(z)}{\partial w_j\partial w_k}-\frac{\partial^2 F_{w_2}(z)}{\partial w_j\partial w_k}\Vert\leq \Vert w_1 -w_2\Vert m_3(z).
\]
So, there exists an integrable function $g\left((Y_1,Z_1),\cdots,(Y_n,Z_n)\right)$ such that, for all $w_1$ and $w_2$ in $\cal W$
\[
\Vert HU_n(w_1)- HU_n(_2)\Vert\stackrel{a.s.}{\leq} \Vert w_1 -w_2\Vert g\left((Y_1,Z_1),\cdots,(Y_n,Z_n)\right).
\]
It follows from this inequality that
\[
\int_0^1HU_n(\hat W_n+u(\hat W_n-w^0)du-HU_n(w^0)\stackrel{a.s.}{\rightarrow}0.
\]
Finally, Theorem 2 is an obvious consequence of this last equation and lemma 2. $\blacksquare$ 

\vskip 3mm

\noindent BIBLIOGRAPHY

\vskip 3mm

\noindent Comrie, A.C. (1997).  Comparing neural networks and regression models for ozone forecasting. {\it Air and Waste Management Association}, {\bf 47}, 653--663.

\vskip 3mm

\noindent Cottrell, M., Girard, B., Girard, Y., Mangeas, M. and Muller, C.  (1995).
Neural modelling for time series: a statistical stepwise method for weight elimination.
{\it IEEE Transaction on Neural Networks}, {\bf 6}, 1355--1364.

\vskip 3mm

\noindent Duflo, M. (1997)
{\it Random iterative models}. Berlin: Springer-verlag.

\vskip 3mm

\noindent Gallant, R. A. (1987).
{\it Non linear statistical models}.  New York: J. Wiley and Sons.

\vskip 3mm

\noindent Gardner, M. W. and Dorling, S. R. (1998).
Artificial neural networks, the multilayer Perceptron. A review of applications in the atmospheric sciences. {\it Atmospheric Environment}, {\bf 32:14/15},  2627--2636.

\vskip 3mm

\noindent Gourieroux, C., Monfort, A. and Trognon, A. (1984). 
Pseudo maximum likelihood methods: Theory.
{\it Econometrica}, {\bf 52:3}, 681--700. 

\vskip 3mm

\noindent Haykin, S. (1999)
{\it Neural networks: A comprehensive foundation}. New Jersey: Prentice Hall.

\vskip 3mm

\noindent Ljung, L. (1999)
{\it System identification: Theory for the user}. New Jersey: Prentice Hall.

\vskip 3mm

\noindent L\"utkepohl, H. (1993)
{\it Introduction to Multiple Time Series Analysis}.
Berlin: Springer-Verlag.

\vskip 3mm

\noindent Magnus, J. and Neudecker, H. (1988).
{\it Matrix differential calculus with applications in statistics and econometrics}.
New York: J. Wiley and Sons. 

\vskip 3mm

\noindent Press, W. H., Teukolsky, S. A., Vetterling, W. T. and Flannery, B. P. (1992).
{\it Numerical recipes in C: The art of scientific computing}.
Cambridge: Cambridge University Press.

\vskip 3mm

\noindent Rumelhart, D. E.,  Hinton, G. E.  and Williams, R. J. (1986).
Learning internal representations by error propagation. In:  D.E. Rumelhart, J.L. McClelland and PDP Research Group, ed.,
{\it Parallel distributed processing Vol. 1}. MIT Press. 

\vskip 3mm

\noindent Rynkiewicz, J. (2003).
Estimation of Multidimensional Regression Model with Multilayer Perceptron.  In J. Mira and A. Prieto, ed. {\it Lecture Notes in Computer Science,  proc. IWANN'2003, Vol. 1} Springer., 310--317. 

\vskip 3mm

\noindent Sussmann, H. J. (1992).
Uniqueness of the weights for minimal feedforward nets with a given input-output Map.
{\it Neural Networks}, {\bf 5},  589--593.

\vskip 3mm

\noindent van der Vaart, A. W. (1998).
{\it Asymptotic statistics}.
Cambridge: Cambridge University Press.

\vskip 3mm

\noindent White, H. (1992).
{\it Artificial neural networks}. Cambridge: Blackwell.

\vskip 3mm

\noindent Yao, J. F. (2000).
On least squares estimation for stable nonlinear AR processes.
{\it The Annals of the Institute of Mathematical Statistics}, {\bf 52}, 316--331.

\end{document}